\documentclass[10pt,reqno]{amsart}

\usepackage{amsmath}
\usepackage{amsthm}
\usepackage{amsfonts}
\usepackage{amssymb,latexsym}
\usepackage[all]{xy}
\usepackage{graphicx}
\usepackage[mathscr]{eucal}
\usepackage{verbatim}
\usepackage{hyperref}
\makeatletter
\def\HyPsd@CatcodeWarning#1{}
\makeatother

\theoremstyle{plain}
    \newtheorem{thm}{Theorem}[section]
    \newtheorem{prop}[thm]{Proposition}
    \newtheorem{lemma}[thm]{Lemma}

\theoremstyle{definition}

\theoremstyle{remark}
    \newtheorem{rem}[thm]{Remark}

\newcommand{\rar}{\ensuremath{\rightarrow}}
\newcommand{\lrar}{\ensuremath{\longrightarrow}}

\newcommand{\la}{\langle}
\newcommand{\ra}{\rangle}
\newcommand{\Hom}{\textup{Hom}}
\newcommand{\StMod}{\textup{StMod}}
\newcommand{\stmod}{\textup{stmod}}

\newcommand{\uHom}{\underline{\Hom}}

\newcommand{\lstk}[1]{\stackrel{#1}{\longrightarrow}}

\newcommand{\Tate}{\widehat{\textup{H}}}

\newcommand{\thick}{\textup{thick}}

\newcommand{\Ind}{\textup{Ind}}

\newcommand{\uInd}{\underline{\Ind}}

\newcommand{\rad}{\text{rad}}
\newcommand{\soc}{\text{soc}}

\newcommand{\loc}{\textup{loc}}

\begin{document}

\title{Freyd's generating hypothesis for groups with periodic cohomology}
\date{\today}

\author{Sunil K. Chebolu}
\address{Department of Mathematics \\
Illinois State University \\
Normal, IL 61790 USA} \email{schebol@ilstu.edu}

\author{J. Daniel Christensen}
\address{Department of Mathematics \\
University of Western Ontario \\
London, ON N6A 5B7, Canada}
\email{jdc@uwo.ca}

\author{J\'{a}n Min\'{a}\v{c}}
\address{Department of Mathematics\\
University of Western Ontario\\
London, ON N6A 5B7, Canada}
\email{minac@uwo.ca}


\begin{abstract}
Let $G$ be a finite group and let $k$ be a field whose characteristic $p$ divides
the order of $G$.
Freyd's generating hypothesis for the stable module category of
$G$ is the statement that a map between finite-dimensional
$kG$-modules in the thick subcategory generated by $k$ factors through a
projective if the induced map on Tate cohomology is trivial. We show that if $G$
has periodic cohomology then the generating hypothesis holds if and only if the Sylow
$p$-subgroup of $G$ is $C_2$ or $C_3$. We also give some other conditions that are equivalent to the GH
for groups with periodic cohomology.
\end{abstract}

\maketitle
\thispagestyle{empty}

\tableofcontents

\newpage

\section{Introduction}

Motivated by the celebrated generating hypothesis (GH) of Peter Freyd in homotopy theory \cite{freydGH} and its
analogue in the derived category of a commutative ring \cite{GH-D(R), keir}, we have formulated in \cite{CCM2}
the analogue of Freyd's GH in the stable module category $\stmod(kG)$ of a finite
$p$-group $G$, where $k$ is a field of characteristic $p$.
(The stable module category is the tensor triangulated category obtained from the category of finitely
generated left $kG$-modules by killing the projective modules.)
In this setting, the GH is the statement that any map that induces the trivial
map in Tate cohomology is trivial in the stable module category $\stmod(kG)$
(i.e., factors through a projective).
In \cite{CCM3} we showed that the only non-trivial $p$-groups for which is this true
are $C_2$ and $C_3$.
The goal of the current project is to describe the analogue of this hypothesis for arbitrary finite groups
and determine for which groups it is true.
It turns out that the above formulation of the GH is not appropriate for arbitrary finite groups, for,
in general, a finite group $G$ can admit a non-projective $kG$-module whose Tate cohomology is trivial.
Clearly the identity map on such a module will disprove the GH, so
it is unreasonable to expect Tate cohomology to detect all non-trivial maps in $\stmod(kG)$.
As we justify in Section~\ref{identityghost},
instead one has to restrict to the thick subcategory $\thick_G(k)$ generated by $k$ in $\stmod(kG)$. (This
is the smallest full subcategory of $\stmod(kG)$ that contains $k$ and closed under exact triangles and
direct summands.)  So the modified  GH for a group ring $kG$ is the statement that Tate cohomology detects
all non-trivial maps in $\thick_G(k)$, i.e.\ that the Tate cohomology functor
\begin{eqnarray*}
 \thick_G(k) & \lrar & \Tate^*(G,k)\text{-modules}  \\
   M & \longmapsto & \Tate^*(G, M)
\end{eqnarray*}
is faithful. If $G$ is a $p$-group, there is only one simple $kG$-module, namely the trivial module $k$,
consequently $\thick_G(k) = \stmod(kG)$.  Therefore this modified GH
agrees with the aforementioned version of the GH for $p$-groups.
In this paper we  determine those finite groups with periodic cohomology for which the modified GH holds.
Our  results can be summarised in:


\begin{thm} \label{maintheorem}
Let $G$ be a non-trivial
finite group that has periodic cohomology and let $k$ be a field of characteristic $p$ that divides the
order of $G$. Then the following are equivalent.
\begin{enumerate}
\item  The Sylow $p$-subgroup of $G$ is either $C_2$ or $C_3$.
\item  The Tate cohomology functor detects all non-trivial maps in $\thick_G(k)$. That is, the GH  holds for $kG$.
\item  Every module in $\thick_G(k)$ is a direct sum of suspensions of $k$.
\item  The Tate cohomology functor detects all non-trivial maps in the stable category $\StMod(B_0)$ of all modules in
the principal block $B_0$ of $kG$.
\item  Every module in $\StMod(B_0)$ is a direct sum of suspensions of $k$.
\end{enumerate}
\end{thm}

It follows that we can make equivalent statements for any full subcategory
which lies between $\thick_{G}(k)$ and $\StMod(B_0)$, such as $\stmod(B_0)$
and $\loc_G(k)$, the localizing subcategory generated by $k$.
It also follows that $\thick_{G}(k) = \stmod(B_{0})$ and $\loc_{G}(k) = \StMod(B_0)$.

Maps of $kG$-modules that induce the trivial map in Tate cohomology are called \emph{ghosts}. In this terminology,
our main result (the equivalence $(1) \iff (2)$ of the above theorem)
states that there are no non-trivial ghosts in $\thick_G(k)$ if and only if the Sylow $p$-subgroup is
$C_2$ or $C_3$.

It is worth pointing out that the GH for $kG$ depends only on $G$ and the characteristic of $k$.
This is clear from the equivalence $(1) \iff (2)$, but is not \emph{a priori} obvious.

Although we have generalised our result for $p$-groups from~\cite{CCM3}, we should stress that our proof in
\cite{CCM3} does not directly generalise.  Several obstacles and subtle issues that arise
in studying the GH for non-$p$-groups are illustrated in Section~\ref{sec:examples} where we work out some
examples of the GH in detail. One new additional technique used here is block theory. In particular, we make good use of
the main theorems of Brauer and the Green correspondence along with some knowledge of the structure of modules
in the principal block for groups with a cyclic normal Sylow $p$-subgroup via Brauer trees.

In work with Carlson~\cite{CarCheMin2} the first and third authors have disproved the GH for
groups with non-periodic cohomology using techniques from Auslander-Reiten theory and support varieties,
and have thus extended all results in this paper to cover the general case, i.e., without any restrictions
on the finite group $G$. Combined with the results of this paper, this gives a complete classification
of the group algebras of finite groups for which the GH holds.
Some related questions which are motivated by the GH have also been studied in~\cite{CCM}.

The paper is organised as follows.
In Section~\ref{sec:preliminaries} we recall several results from
representation theory which are used in the later sections.  We also
prove that (2) and (3) above are equivalent.
Section~\ref{sec:examples} contains a few important examples which
illustrate some issues that arise when studying the GH for
non-$p$-groups.
The main steps in the proof of Theorem~\ref{maintheorem} occupy
Sections~\ref{sec:ifpart} and~\ref{sec:onlyifpart}. In
Section~\ref{sec:ifpart} we show that (1) implies (3) and
in Section~\ref{sec:onlyifpart} we show that (3) implies (1).
The equivalence of (4) and (5) with the other statements is shown in
Section~\ref{sec:preliminaries}.
The reader who is only interested in the proof of the main theorem may
skip Sections~\ref{sec:preliminaries} and~\ref{sec:examples}, referring
to Section~\ref{sec:preliminaries} when necessary.

All groups in this paper are non-trivial finite groups and the characteristic $p$ of the field always divides
the order of $G$. We work in the stable module category of $kG$ and we freely use standard facts about this category
which can be found in \cite{carlson-modulesandgroupalgebras}.

\subsection*{Acknowledgements} The first and third authors carried out some of this work at the PIMS algebra summer school
at the University of Alberta.  They would like to thank  the University of Alberta,
the organisers (A.~Adem, J.~Kuttler and A.~Pianzola)  of the  summer school, and V.~Chernousov for their hospitality.
We are grateful to J.~Carlson for his work in \cite{CarCheMin} which inspired us to consider almost split sequences
in the last section.
Calculations using Peter Webb's \texttt{reps} package~\cite{reps} for GAP~\cite{GAP}
were useful in our search for non-trivial ghosts.

\section{Some results from representation theory} \label{sec:preliminaries}

In this section we collect some known results from representation theory which we will need in the
sequel.

\subsection{Periodic cohomology}

We say that $kG$, or simply $G$ when there is no confusion, has \emph{periodic cohomology} if there is a positive
integer $d$ such that $\Omega^d k$ is stably isomorphic to $k$.  When
this is the case, the \emph{period} is the smallest such $d$.
It is a well-known fact due to E.~Artin and Tate~\cite[p.\ 262]{CarEil} that a finite group $G$ has
periodic cohomology over a field $k$ of characteristic $p$ if and only if the Sylow $p$-subgroup of $G$ is cyclic or a
generalised quaternion group.

We begin with a proposition which forms the backbone of our analysis.

\begin{prop}[\cite{CCM2}] \label{prop:crucial}
Let $G$ be a finite group with periodic cohomology.
Then the GH holds for $kG$ if and only if every
module in $\thick_G(k)$ is a sum of suspensions of $k$.  In particular, the GH holds for $kG$
if and only if every indecomposable non-projective
$kG$-module in $\thick_G(k)$ is stably isomorphic to $\Omega^i k$ for some $i$.
\end{prop}

\begin{proof} We sketch a proof here; more details can be found in~\cite{CCM2}.  Let $M$ be in
$\thick_G(k)$.
Since the trivial representation is periodic,
a ghost out of $M$ can be constructed in $\thick_G(k)$ using a triangle of the form
\[ \underset{\text{finite sum}}{\bigoplus} \Omega^i k \lrar M \lstk{f} U_M. \]
If the GH holds for $kG$, then $f$ must vanish.
Thus the above triangle splits, and so
$M$ is a retract of $\oplus \, \Omega^i k$. Since $M$ is finite-dimensional, it follows
from the Krull-Schmidt theorem that $M$ is a sum of suspensions of $k$.
The converse is immediate.
\end{proof}

Thus the GH holds if and only if the number of indecomposable non-projective
$kG$-modules in $\thick_G(k)$ is equal to the period.
The next two results give us tools for computing these quantities.

\begin{thm}[Swan~\cite{Swan-period}]\label{thm:Swan}
Let $G$ be a finite group with periodic cohomology.
When $p = 2$, the period is $1$, $2$ or $4$ when the Sylow $2$-subgroup is $C_2$, $C_{2^r}$
($2^r > 2$), or $Q_{2^n}$, respectively. When $p$ is odd and the Sylow $p$-subgroup is $C_{p^r}$, the period is $2 \Phi_p$,
where $\Phi_p$ is the number of automorphisms of $C_{p^r}$ that are given by conjugation by elements in $G$.
\end{thm}

\begin{thm}\label{thm:CR}
Let $G$ be a finite group with cyclic Sylow $p$-subgroup of order $p^r$,
and let $s$ be the number of simple $kG$-modules.
Then the number of indecomposable non-projective $kG$-modules is $s(p^r-1)$.
Moreover, if $B$ is a block of $kG$
and $e$ is the number of simple modules lying in $B$,
then the number of indecomposable non-projective $kG$-modules lying in $B$ is $e(p^r-1)$.
\end{thm}

\begin{proof}
The first statement is a simplified version of~\cite[Prop.~20.11]{Cur-Rei-I}.
The second statement follows from the detailed structure given there, using the fact that a module $M$
lies in the block $B$ if and only if each composition factor of $M$ lies in $B$.
\end{proof}

While Theorem~\ref{thm:CR} doesn't deal directly with $\thick_G(k)$, we
will use it in Section~\ref{sec:ifpart} to show that every $kG$-module in the principal
block is a sum of suspensions of $k$.

\subsection{Partial proof of Theorem~\ref{maintheorem}}

We are now ready to prove all equivalences of Theorem~\ref{maintheorem},
using the results of Sections~\ref{sec:ifpart} and~\ref{sec:onlyifpart}.
The implications $(5) \implies (4) \implies (2)$ are clear.
In Proposition~\ref{prop:crucial} we have seen that $(2) \iff (3)$.
In Section~\ref{sec:onlyifpart} we prove $(3) \implies (1)$.
Thus it remains to prove $(1) \implies (5)$.  In Section~\ref{sec:ifpart}
we show that $(1) \implies (3)$, but we in fact show a stronger result:
$(1)$ implies that every module in $\stmod(B_0)$ is a sum of suspensions of $k$.
Now a result of Ringel and Tachikawa~\cite{ringel-tachikawa} states that
if $G$ has finite representation type (i.e., the Sylow $p$-subgroups are
cyclic), then every $kG$-module is a direct sum of finite-dimensional
$kG$-modules.
It follows that when $(1)$ holds, every module in $\StMod(B_0)$ is a
sum of modules in $\stmod(B_0)$, and so $(5)$ follows.
\qed

The only places in this paper where we used the assumption that $G$
has periodic cohomology are in ruling out the possibility that the
Sylow $p$-subgroup is a dihedral $2$-group and in Proposition~\ref{prop:crucial}.
(We use periodicity in Section~\ref{sec:ifpart}, but there it follows from
the assumption that the Sylow $p$-subgroup is $C_{2}$ or $C_{3}$.)
Thus we can make the following statement, without the hypothesis that
$G$ has periodic cohomology.

\begin{thm}
Let $G$ be a group whose Sylow $p$-subgroup is not a dihedral $2$-group.
Then the Sylow $p$-subgroup of $G$ is $C_{3}$ if and only if
every module in $\thick_{G}(k)$ is a sum of suspensions of $k$.
\end{thm}

Of course, if $p$ is odd, then the first condition on the Sylow $p$-subgroup
can be omitted. The case $p = 2$ is completed in~\cite{CarCheMin2}.

\section{Examples} \label{sec:examples}
In this section we discuss some  examples which will help the reader get some insight into the GH.

\subsection{Non-trivial identity ghosts} \label{identityghost}

It is well-known (see, e.g.,~\cite{keir}) that the right setting for
the GH in a general triangulated category is the thick subcategory
generated by the distinguished object (in our case, the trivial
representation $k$).
For the stable module category of a group algebra, it is
difficult  to illustrate why this is the right choice, since our
main result implies that the GH holds in $\thick_G(k)$ if and only if
it holds in any full subcategory containing $\thick_G(k)$ and
contained in $\StMod(B_0)$, where $B_0$ is the principal block.
Moreover, when the GH holds, we show that $\thick_G(k) = \stmod(B_0)$.
However, we can study identity maps which are ghosts in order to
get some insight into this issue.

The key point is that, in general, there can be non-projective modules with trivial Tate cohomology.
Clearly the identity map on such a module will be a non-trivial ghost.  Examples of such modules abound.
For instance, if there is a non-projective indecomposable module $M$
that does not belong to the principal block $B_0$, then this gives an example.
So clearly one needs to restrict to the principal block.
Moreover, if $\thick_{G}(k)$ is a proper subcategory of $\stmod(B_0)$,
then work of Benson, Carlson and Robinson~\cite{Benson-NY, BCR-90} shows that
there is an indecomposable non-projective module which is in $\stmod(B_0)$ but outside of $\thick_{G}(k)$
and has trivial Tate cohomology.

In contrast, we show that there are no non-trivial identity ghosts
in the thick subcategory generated by $k$.
This gives some evidence that $\thick_G(k)$ is the ``right'' category
in which to study the GH.

\begin{prop} \label{prop:identityghostinthick(k)}
Let $M$ be in $\thick_G(k)$.  If the identity map
$M \rar M$ is a ghost, then it is trivial in $\stmod(kG)$.
\end{prop}

\begin{proof}
This is a standard thick subcategory argument. Consider the full subcategory of all
modules $X$ in $\stmod(kG)$ which have the property that
$\uHom(\Omega^i X, M) = 0$ for all integers $i$. It is
straightforward to verify that this subcategory is closed under retractions
and exact triangles. It contains the trivial representation by
hypothesis.  Thus it contains the thick subcategory generated by $k$,
and hence contains $M$. In particular the identity map on $M$ is trivial.
\end{proof}

In some favourable cases, even when $G$ is not a $p$-group, $\thick_G(k)$ can be the whole of $\stmod(kG)$.
The GH for such groups can be easily attacked using the restriction-induction technique of~\cite{CCM3}. We illustrate this in
the example of $A_4$.

\subsection{The alternating group $A_4$ when $p =2$} \label{A_4}
Let $k$ be a field of characteristic $2$ and
consider the alternating group $A_4$. This is a group of
order $12$ and is generated by $x$, $y$ and $z$ which satisfy the
relations $x^2 =  y^2 = (xy)^2 = 1 = z^3, zxz^{-1} = y$ and $zyz^{-1} = xy$.
Using these relations, one can show that the centraliser of every element of order $2$ is $2$-nilpotent.
Work of Benson, Carlson and Robinson~\cite{Benson-NY, BCR-90} then implies
that $\thick_G(k) = \stmod(B_0)$.
Moreover the principal idempotent can be shown to be $1$, so we in fact have
$\thick_{A_4}(k) = \stmod(kA_4)$.

%

Now the subgroup of $A_4$ generated by $x$ and $y$ is the Klein
four group $V_4$.
%
So the Sylow $2$-subgroup is $V_4$.
By \cite{CCM3}, we know that the GH fails for $V_4$. So the induction of a non-trivial ghost over $kV_4$
will give a non-trivial ghost (see \cite[Prop.~2.1]{CCM3}) over $kA_4$, thus disproving the GH for $kA_4$.

\begin{rem}
The induction  functor $\uInd\colon \stmod(kH) \rar \stmod(kG)$ does not in general send $\thick_H(k)$ into $\thick_G(k)$.
For example, if $\mathbb{F}_3$ is the trivial $\mathbb{F}_{3}C_3$ module, then it can be shown that the induced
$\mathbb{F}_{3}(C_2 \times C_3)$-module $\mathbb{F}_3{\uparrow^{C_2 \times C_3}}$ does not belong to the thick
subcategory $\thick_{C_2 \times C_3}(\mathbb{F}_3)$. Since the right domain for the GH is $\thick_G(k)$,
the above induction strategy does not generalise to arbitrary finite groups.
\end{rem}

\subsection{The symmetric group $S_3$ when $p = 3$} \label{S_3}

In this section we prove that the GH holds in $\thick_{S_3}(k)$ when $k$
has characteristic $3$.  The argument we give here is a model for
the general argument we give in Section~\ref{sec:ifpart}, and
also illustrates Theorem~\ref{thm:CR}.

The group $S_3$ has presentation
$\la x, y \, | \, x^3 = 1 = y^2, yxy^{-1} = x^{-1} \ra$.
Define elements $e_1 = (1 - y)/2$ and $e_2 = (1 + y)/2$ in $A = kS_3$.
Then $e_1 + e_2 = 1$ and it is
a straightforward exercise to show that $e_1$ and $e_2$ are orthogonal idempotents in $A$, i.e.\
$e_1^2 = e_1$, $e_2^2 = e_2$ and $e_1 e_2 = 0 = e_2 e_1$.  The principal indecomposable
modules $Ae_1$ and $Ae_2$ (both $3$-dimensional) have composition series of length 3:
\[ Ae_1 \supsetneq A(x-1)e_1 \supsetneq A(x-1)^2e_1 \supsetneq 0 \]
\[ Ae_2 \supsetneq A(x-1)e_2 \supsetneq A(x-1)^2e_2 \supsetneq 0 . \]
These six modules form a complete set of indecomposable $kS_3$-modules; see \cite[\S~64]{CR}.
Moreover, $Ae_1$ and $Ae_2$ are the indecomposable projectives over the simple modules $A(x-1)^2e_1$
and $A(x-1)^2e_2$ respectively. The structure of the simples is as follows:
$A(x-1)^2e_2 = k$, the trivial representation, and $A(x-1)^2e_1 = k_{-1}$, on which $x$ acts trivially and
$y$ by multiplication by $-1$. We now leave it as an amusing exercise for the reader to show that
\begin{align*}
k &\cong A(x-1)^2e_2, \\
\Omega k &\cong A(x-1)\,\,e_2, \\
\Omega^2 k &\cong A(x-1)^2e_1 \, (= k_{-1}), \\
\Omega^3 k &\cong A(x-1)\,\,e_1, \ \  \text{and }\\
\Omega^4 k &\cong k.
\end{align*}
So $k$ has period $4$, which agrees with the answer we get from
Swan's formula (Theorem~\ref{thm:Swan}): $2\Phi_3 = 2(2)  = 4$.
This also shows that every indecomposable non-projective $kG$-module is isomorphic to
$\Omega^i k$ for some $i$, and so the GH holds for $kS_3$.

This example suggests that the GH for non-$p$-groups is both subtle and interesting.

\section{Groups with periodic cohomology for which the GH holds} \label{sec:ifpart}

In this section we show that if the Sylow $p$-subgroup of $G$ is either $C_2$
or $C_3$, then every module in $\stmod(B_{0})$ is a sum of suspensions of $k$,
where $B_{0}$ is the principal block of $kG$.
 From this it follows that the GH holds for $kG$.

We next give some results which will be used in the proof.

\subsection{Field extensions}

\begin{lemma}
Let $L$ be an extension of $k$ and let $G$ be a finite group.
Then the principal block of $LG$ is $L \otimes_k B_{0}$, where
$B_0$ is the principal block of $kG$.
Moreover, if every module in $\stmod(L \otimes_k B_{0})$ is a
sum of suspensions of $L$, then every module in $\stmod(B_{0})$
is a sum of suspensions of $k$.
\end{lemma}

\begin{proof}
The statement about the principal block of $LG$ follows from the
fact that the principal idempotent depends only on the characteristic
of the field (see, e.g.,~\cite{kuelshammer}).

To prove the second statement, note that the functor
\[
  L \otimes_k - : \StMod(kG) \longrightarrow \StMod(LG)
\]
is faithful, triangulated and sends ghosts to ghosts.  It restricts to a functor
\[
  L \otimes_k - : \StMod(B_{0}) \longrightarrow \StMod(L \otimes_{k} B_{0}).
\]
Let $M$ be a $kG$-module in $\stmod(B_{0})$.  Consider the triangle
\[
  \bigoplus_{i\in\mathbb{Z}}\,\bigoplus_{\eta\in\uHom(\Omega^{i}k, M)} \Omega^i k \lrar M \lstk{\Phi_M} U_M
\]
in $\StMod(B_{0})$.  If every $LG$-module in $\stmod(L \otimes_{k} B_{0})$
splits as a sum of suspensions of $L$, then $L \otimes_{k} \Phi_{M}$ is
stably trivial, and so $\Phi_{M}$ is stably trivial.
Thus, using Krull-Schmidt, $M$ splits as a sum of suspensions of $k$.
\end{proof}

Thus we can assume that $k$ is algebraically closed, and we do so for the remainder of this section.
This is convenient because we cite~\cite{Alperin-book} in
Sections~\ref{subsec:normal} and~\ref{subsec:C_3}, and that
reference makes the assumption that $k$ is algebraically closed.

\subsection{Direct products}

\begin{lemma} \label{lemma:PxP'}
Let $G$ be a finite group that is a product of two groups: $G = A \times B$. Assume that $p$ does not
divide the order of $B$.  Then the restriction functors
\begin{align*}
 \stmod(kG)  &\to \stmod(kA)
\intertext{and}
 \StMod(kG)  &\to \StMod(kA)
\end{align*}
are tensor triangulated equivalences of categories.
\end{lemma}

This lemma is well-known, but we give a proof here for the reader's convenience.

\begin{proof}
The restriction functors are easily seen to be tensor triangulated functors.
That is, they preserve suspension, cofibre sequences and tensor products,
and they send the unit object $k$ to the unit object $k$.
Since any $kA$-module can be viewed as a $kG$-module with a trivial
action of $B$, the restriction functors are full and essentially surjective.
We only need to show that they are faithful.
This is true for any subgroup $A$ whose index in $G$ is invertible in $k$,
since the composite of the restriction map
\def\downA{{\downarrow}_{A}}
\[
   \uHom_{G}(M, N) \longrightarrow \uHom_{A}(M \downA, N \downA)
\]
with the transfer map
\[
  \uHom_{A}(M \downA, N \downA) \longrightarrow \uHom_{G}(M, N)
\]
is multiplication by $[G : A]$.
\end{proof}

It follows that $\thick_{G}(k)$ is equivalent to $\thick_{A}(k)$,
and one can also show that $kG$ and $kA$ have isomorphic principal
blocks.

\begin{rem}\label{rem:S_3}
This result cannot be generalised to semi-direct products.
The example to keep in mind is $kS_3 = k(C_{3} \rtimes C_{2})$, where the characteristic
of $k$ is $3$. By Swan's formula~(Theorem~\ref{thm:Swan}) or the computations in Section~\ref{S_3},
the trivial representation $k$ has period $4$ in $\thick_{S_3}(k)$
and has period $2$ in $\thick_{C_3}(k)$. In particular,
\[ \thick_{S_3}(k) \not \cong \thick_{C_3}(k). \]
So while the point of this paper is to show that the GH is determined
by the Sylow $p$-subgroup, it is not because the relevant thick subcategories
are equivalent.
\end{rem}

\subsection{Reduction to the normal case}\label{subsec:normal}

We now use results from block theory to
show that when the Sylow $p$-subgroup $D$ of $G$ is $C_p$,
we can reduce to the case where $D$ is normal.
The relevant background material can be found in~\cite{Alperin-book,ben-1},
for example.

\begin{thm}\label{thm:reductiontonormalcase}
Let $G$ be a group which has a cyclic Sylow $p$-subgroup $D$,
let $D_1$ be the unique subgroup of $D$ that is isomorphic to $C_p$
and let $N_1 = N_G(D_1)$.
Then there is a tensor triangulated equivalence of categories
\[ \stmod(B_0) \cong \stmod(b_0), \]
where $B_0$ is the principal block of $kG$ and $b_0$ the principal block of $kN_1$.
\end{thm}

When $D$ is $C_p$, then $D_1 = D$ and so $D$ is also the Sylow $p$-subgroup
of $N_1$ and is normal in $N_1$.

\begin{proof}
Recall that $D$ is the defect group of the principal block.
Since $D C_G(D) = C_G(D) \le N_1$,
Brauer's third main theorem says that the block $b_0^G$ corresponding to
the principal block $b_0$ of $kN_1$ is the principal block $B_0$ of $kG$.
So by~\cite[pp.~124--125]{Alperin-book}, there is an equivalence of categories
\[ \stmod(B_0) \cong \stmod(b_0). \qedhere\]
\end{proof}

By Theorem~\ref{thm:reductiontonormalcase},
we know that if the Sylow $p$-subgroup $H$ of $G$ is isomorphic to $C_p$,
then the stable categories of the principal blocks of $kG$ and $kN_G(H)$ are equivalent.
So we can assume without loss of generality that $H$ is normal in $G$.

\subsection{The Sylow $p$-subgroup is $C_2$} If $H = C_2$ is normal in $G$, then it is actually central in $G$.
By the Schur-Zassenhaus Theorem it follows that $G = C_2 \times L$ for some group $L$ which has odd order.  Then by
Lemma~\ref{lemma:PxP'} we have that $\stmod(kG)$ is equivalent to $\stmod(kC_2)$ as tensor triangulated
categories.
By the main result of~\cite{CCM3}, every module in $\stmod(kC_{2})$ is
a sum of suspensions of $k$, so the same is true in $\stmod(kG)$.
In particular, this is true for the principal block.

\subsection{The Sylow $p$-subgroup is $C_3$}\label{subsec:C_3}
Let  $H = C_3$ be normal in $G$. Now consider the map
\begin{eqnarray*}
\Xi \colon G & \lrar & \text{Aut}(C_3) \cong C_2 \\
     g & \longmapsto & g(-)g^{-1}
\end{eqnarray*}
There are only two possibilities for the image of $\Xi$:

\subsubsection*{Case 1:}
The image of $\Xi$ is trivial. In this case, exactly as before, $G = C_3 \times L$ for some group
$L$ whose order is not divisible by $3$. So by Lemma~\ref{lemma:PxP'} we have that $\stmod(kG)$
is equivalent to $\stmod(kC_3)$ as tensor triangulated categories.
By the main result of~\cite{CCM3}, every module in $\stmod(kC_{3})$ is
a sum of suspensions of $k$, so the same is true in $\stmod(kG)$.
In particular, this is true for the principal block.

\subsubsection*{Case 2:}
The image of $\Xi$ is $C_2$. Then the centraliser $C_G(C_3)$ has index $2$ in $G$. In this case,
$\Phi_3$, the number of automorphisms of $C_3$ given by conjugation by elements of $G$, is equal to $2$.
By Theorem~\ref{thm:Swan}, $k$ has period $2\Phi_3 = 4$.
Thus it is enough to show that there are exactly four
indecomposable non-projective $kG$-modules in the principal block.
By Theorem~\ref{thm:CR}, we know that the number of indecomposable non-projective $kG$-modules in the principal
block is twice the number of simple $kG$-modules in the principal block.  Combining these two, we just have to show that
there are only two simple $kG$-modules in the principal block.

Let $P$ be the indecomposable projective module over $k$, that is $P/\rad(P) \cong k$.
Let $W$ be the module $\rad(P)/\rad^2(P)$. Then the set of all simple $kG$-modules
in the principal block is
\[\{k,\, W,\, W \otimes W,\, W \otimes W \otimes W,\, \cdots \}.\]
This fact can be found in~\cite[Exercise~13.3]{Alperin-book}, for instance.
We will be done if we can show that $k \ncong W$ and
$W \otimes W = k$ because then we will have exactly two
simple $kG$-modules in the principal block, namely $k$ and $W$.  These two facts will become clear once we give the
explicit structure  of $W$. Write $G = C_3 L$, where $L$ is a complement of $C_3$, which exists by the Schur-Zassenhaus theorem,
and let $x$ be a generator of $C_3$. It can be shown \cite[p.~37]{Alperin-book} that $W$
is a one-dimensional module generated by $v$ such that $x(v) = v$, and for $h$ in $L$, $h(v) = v$ if $h$
belongs to $C_G(C_3)$ and  $-v$ if $h$ does not belong to $C_G(C_3)$. Since $C_G(C_3)$ has index $2$, there are
elements outside $C_G(C_3)$ which do not fix $v$, and therefore $W$ is not isomorphic to $k$.  The fact that
$W \otimes W \cong k$ is clear since
\begin{align*}
x(v \otimes v) &= x v \otimes x v  =    v \otimes v   \\
h(v \otimes v) &= h v \otimes h v  =  \pm v \otimes \pm v  = v \otimes v .
\end{align*}
This shows that there are exactly two simple $kG$-modules in the principal block. So we are done.

\section{Groups with periodic cohomology for which the GH fails}  \label{sec:onlyifpart}

In this section we show that for a group $G$ which has periodic cohomology, the GH fails  whenever the Sylow
$p$-subgroup of $G$ is not $C_2$ or $C_3$.  In view of Proposition~\ref{prop:crucial}, in order to
disprove the GH for these groups  we have to  show that there is a module in $\thick_G(k)$
that is not stably isomorphic to a direct sum of suspensions of $k$.
We will show that the middle term of an almost split sequence has this property.

We recall the standard almost split sequence for the reader. Let $G$ be any finite group and
let $P$ be the indecomposable projective module over $k$, that is, $P/\rad{P} \cong k$.
Since $kG$ is a symmetric algebra, we also have $\soc{P} \cong k$. The quotient $\rad{P}/\soc{P}$ is
called the heart $H_G$ of $G$. It occurs as a summand in the middle term of the standard almost split sequence
\[ 0 \lrar \rad{P} \lrar  H_G \oplus P \rar P/ \soc{P} \lrar 0. \]
This sequence can also be written as
\begin{equation}\label{eq:ses}
  0 \lrar \Omega^1 k \lrar  H_G \oplus P \rar \Omega^{-1} k \lrar 0 .
\end{equation}
It is a non-trivial result of Webb \cite[Thm.~E]{Webb} that $H_G$ is an indecomposable $kG$-module provided the Sylow
$p$-subgroup of $G$ is not a dihedral $2$-group.
This covers our situation, since for a group with periodic cohomology, the
only dihedral $2$-group that can arise as the Sylow $p$-subgroup is $C_{2}$,
and we are explicitly excluding this possibility.

\begin{thm}\label{thm:heart}
Let $G$ be a group which has periodic cohomology for which the Sylow $p$-subgroup is not $C_2$ or $C_3$.
Then the $kG$-module $H_G$ is an indecomposable non-projective module  in $\thick_G(k)$ that is not
stably isomorphic to $\Omega^i k$ for any $i$. In particular, there is a non-trivial ghost out of
$H_G$ in $\thick_G(k)$, i.e., the GH fails for $kG$.
\end{thm}

\begin{proof}
It is clear from the short exact sequence~\eqref{eq:ses} that $H_G$ belongs to $\thick_G(k)$.
Further,  we know from Webb's theorem stated above  that $H_G$ is an indecomposable $kG$-module.
So we only have to show that $H_G$ is not projective and that it is not stably isomorphic to $\Omega^i k$ for any $i$.
Both of these statements follow easily by comparing dimensions. The key fact to observe is that the dimension of every
projective $kG$-module is divisible by $p^n$,  the order of the Sylow $p$-subgroup of $G$.
(One sees this by restricting the projective module to the  Sylow $p$-subgroup $P$, over which the restriction
becomes a free $kP$-module.) On the other hand, from
the definition of $H_G$, it is clear that  $\dim_k H_G \equiv -2 \mod p^n$.
So if $H_G$ is projective, then $p^n$ should
divide $2$, but that would mean that the Sylow $p$-subgroup is $C_2$, which is a contradiction.
Therefore $H_G$ has to be non-projective.
Using the minimal projective resolution of $k$ and the above fact about dimensions of projective $kG$-modules,
one sees by a straightforward induction on $i$ that $\dim_k \Omega^i k \equiv 1$ or $-1$ mod $p^n$.
If $H_G \cong \Omega^i k$ for some $i$, then it follows that the Sylow $p$-subgroup is either trivial or $C_3$.
Both cases are ruled out by our assumptions.  Therefore $H_G$ is not stably isomorphic to $\Omega^i k$ for any $i$.
The last statement follows from Proposition~\ref{prop:crucial}.
\end{proof}

The last two sections together prove our main theorem that if $G$ has periodic
cohomology, then the GH holds for $kG$ if and only if the Sylow $p$-subgroup of $G$ is either $C_2$ or $C_3$.

\enlargethispage{1pt}


\end{document}